\documentclass[12pt]{amsart}
\usepackage{amsmath,amsthm,amsfonts,amssymb,mathrsfs}
\date{\today}

\usepackage{color}

\usepackage{hyperref}

\newtheorem{theorem}{Theorem}[section]

\newtheorem{proposition}[theorem]{Proposition}
\newtheorem{corollary}[theorem]{Corollary}
\newtheorem{lemma}[theorem]{Lemma}
\theoremstyle{definition}
\newtheorem{remark}[theorem]{Remark}

\begin{document}

\title[On semitopological bicyclic extensions of linearly ordered groups] {On semitopological bicyclic extensions of linearly ordered groups}

\author[Oleg~Gutik and Kateryna Maksymyk]{Oleg~Gutik and Kateryna Maksymyk}
\address{Faculty of Mathematics, National University of Lviv,
Universytetska 1, Lviv, 79000, Ukraine}
\email{o\underline{\hskip5pt}\,gutik@franko.lviv.ua,
ovgutik@yahoo.com, kate.maksymyk15@gmail.com}



\begin{abstract}
For a linearly ordered group $G$ let us define a subset $A\subseteq G$ to be a \emph{shift-set} if for any
$x,y,z\in A$ with $y < x$ we get $x\cdot y^{-1}\cdot z\in A$. We describe the natural partial order and
solutions of equations on the semigroup $\mathscr{B}(A)$ of shifts of positive cones of $A$. We study topologizations of the semigroup $\mathscr{B}(A)$. In particular, we show that for an arbitrary countable linearly ordered group $G$ and a non-empty shift-set $A$ of $G$  every Baire shift-continuous $T_1$-topology $\tau$ on $\mathscr{B}(A)$ is discrete. Also we prove that for an arbitrary linearly non-densely ordered group $G$ and a non-empty shift-set $A$ of $G$, every shift-continuous Hausdorff topology $\tau$ on the semigroup $\mathscr{B}(A)$ is discrete, and hence $\left(\mathscr{B}(A),\tau\right)$ is a discrete subspace of any Hausdorff semitopological semigroup which contains $\mathscr{B}(A)$ as a subsemigroup.
\end{abstract}

\maketitle

\section{Introduction and preliminaries}

We shall follow the terminology of \cite{Carruth-Hildebrant-Koch-1983-1986, Clifford-Preston-1961-1967, Engelking-1989, Fuchs-1963, Haworth-McCoy-1977, Petrich-1984, Ruppert-1984}.

A \emph{semigroup} is a non-empty set with a binary associative
operation. A semigroup $S$ is called \emph{inverse} if for any $x\in
S$ there exists a unique $y\in S$ such that $x\cdot y\cdot x=x$ and
$y\cdot x\cdot y=y$. Such an element $y$ in $S$ is called the
\emph{inverse} of $x$ and is denoted by $x^{-1}$. The map defined on an
inverse semigroup $S$ which maps every element $x$ of $S$ to its
inverse $x^{-1}$ is called the \emph{inversion}.

For a semigroup $S$ by $E(S)$ we denote the set of idempotents in $S$. If $S$ is an inverse semigroup, then
$E(S)$ is closed under multiplication and we shall refer to $E(S)$
as the \emph{band of} $S$.  A \emph{semilattice} is a commutative semigroup of
idempotents.

Let $\mathscr{I}_X$ denote the set of all partial one-to-one
transformations of an infinite set $X$  together with the following
semigroup operation: $x(\alpha\beta)=(x\alpha)\beta$ if
$x\in\operatorname{dom}(\alpha\beta)=\{
y\in\operatorname{dom}\alpha\colon
y\alpha\in\operatorname{dom}\beta\}$,  for
$\alpha,\beta\in\mathscr{I}_X$. The semigroup $\mathscr{I}_X$ is
called the \emph{symmetric inverse semigroup} over the set $X$~(see
\cite{Clifford-Preston-1961-1967}). The symmetric inverse semigroup
was introduced by Wagner~\cite{Wagner-1952} and it plays a major role
in the theory of semigroups.

The bicyclic monoid ${\mathscr{C}}(p,q)$ is the semigroup with
the identity $1$ generated by two elements $p$ and $q$ subjected
only to the condition $pq=1$. The bicyclic monoid is a combinatorial bisimple $F$-inverse semigroup and it plays an important role in the algebraic
theory of semigroups and in the theory of topological semigroups.
For example the well-known O.~Andersen's result~\cite{Andersen-1952}
states that a ($0$--) simple semigroup is completely ($0$--) simple
if and only if it does not contain the bicyclic monoid. The
bicyclic monoid does not embed into stable
semigroups~\cite{Koch-Wallace-1957}.

Recall \cite{Fuchs-1963} that a \emph{partially ordered group}
is a group $(G,\cdot)$ equipped with a translation-invariant partial order  $\leqslant$; in other words, the binary relation $\leqslant$  has the
property that, for all $a, b, g\in G$, if $a\leqslant b$ then
$a\cdot g\leqslant b\cdot g$ and $g\cdot a\leqslant g\cdot b$.

Later by $e$ we denote the identity of a group $G$. The set
$G^+=\{x\in G\colon e\leqslant x\}$ in a partially ordered group $G$
is called the \emph{positive cone}, or the \emph{integral part}, of
$G$ and it satisfies the properties:
\begin{equation*}
  1)~G^+\cdot G^+\subseteq G^+; \quad 2)~G^+\cap (G^+)^{-1}=\{e\}; \quad \hbox{and} \quad 3)~x^{-1}\cdot G^+\cdot x\subseteq G^+ \hbox{~for each~} x\in G.
\end{equation*}
Any subset $P$ of a group $G$ that satisfies the conditions 1)--3)
induces a partial order on $G$ ($x\leqslant y$ if and only if
$x^{-1}\cdot y\in P$) for which $P$ is the positive cone. Elements of the set $G^+\setminus\{e\}$ are called \emph{positive}.

A \emph{linearly ordered} or \emph{totally ordered group} is an
ordered group $G$ whose order relation $\leqslant$ is
total~(see \cite{Birkhoff-1973} and \cite{Clay-Rolfsen-2016}).

From now on we shall assume that $G$ is a non-trivial linearly ordered
group.

For every $g\in G$ the set
\begin{equation*}
    G^+(g)=\{x\in G\colon g\leqslant x\}.
\end{equation*}
The set $G^+(g)$ is called a \emph{positive cone on element} $g$ in
$G$.

For arbitrary elements $g,h\in G$ we consider a partial map
$\alpha_h^g\colon G\rightharpoonup G$ defined by the formula
\begin{equation*}
    (x)\alpha_h^g=x\cdot g^{-1}\cdot h, \qquad \hbox{ for } \; x\in
    G^{+}(g).
\end{equation*}
We observe that Lemma~XIII.1 from \cite{Birkhoff-1973} implies that
for such partial map $\alpha_h^g\colon G\rightharpoonup G$ the
restriction $\alpha_h^g\colon G^+(g)\rightarrow G^+(h)$ is a
bijective map.

We consider the semigroups
\begin{equation*}
    \mathscr{B}(G)=\{\alpha_h^g\colon G\rightharpoonup G\colon g,h\in
    G\} \, \hbox{ and } \,
    \mathscr{B}^+(G)=\{\alpha_h^g\colon G\rightharpoonup G\colon g,h\in
    G^+\},
\end{equation*}
endowed with the operation of the composition of partial maps. Simple verifications show that
\begin{equation}\label{formula-1.1}
\alpha_h^g\cdot \alpha^k_l=\alpha^a_b, \qquad \hbox{ where } \quad
a=(h\vee k)\cdot h^{-1}\cdot g \quad \hbox{ and } \quad b=(h\vee
k)\cdot k^{-1}\cdot l,
\end{equation}
for $g,h,k,l\in G$, and by $h\vee k$ we denote the join of $h$ and $k$ in the linearly ordered set $(G,\leqslant)$. Therefore, property 1) of the positive cone and
condition~\eqref{formula-1.1} imply that $\mathscr{B}(G)$ and
$\mathscr{B}^+(G)$ are subsemigroups of $\mathscr{I}_G$.

By Proposition~1.2 in \cite{Gutik-Pagon-Pavlyk-2011} for a linearly ordered group $G$ the following assertions hold:
\begin{itemize}
  \item[$(i)$] elements $\alpha_h^g$ and $\alpha_g^h$ are inverses
   of each other in $\mathscr{B}(G)$ for all $g,h\in G$
   $($resp., $\mathscr{B}^+(G)$ for all $g,h\in G^+)$;

  \item[$(ii)$] an element $\alpha_h^g$ of the semigroup
   $\mathscr{B}(G)$ $($resp., $\mathscr{B}^+(G))$ is an
   idempotent if and only if $g=h$;

  \item[$(iii)$] $\mathscr{B}(G)$ and $\mathscr{B}^+(G)$ are inverse
   subsemigroups of $\mathscr{I}_G$;

  \item[$(iv)$] the semigroup $\mathscr{B}(G)$ $($resp.,
   $\mathscr{B}^+(G))$ is isomorphic to the set $S_G=G\times G$ $($resp.,
   $S_G^+=G^+\times G^+)$ with the following semigroup operation:
  \begin{equation}\label{formula-1.2}
  (a,b)(c,d)=
  \left\{
    \begin{array}{ll}
      (c\cdot b^{-1}\cdot a,d), & \hbox{if }  \; b<c;\\
      (a,d), & \hbox{if } \; b=c; \\
      (a,b\cdot c^{-1}\cdot d), & \hbox{if } \; b>c,
    \end{array}
  \right.
  \end{equation}
  where $a,b,c,d\in G$ $($resp.,  $a,b,c,d\in G^+$$)$.
\end{itemize}

It is obvious that:
\begin{itemize}
  \item[$(1)$] if $G$ is isomorphic to the additive group of integers
  $(\mathbb{Z},+)$ with usual linear order $\leqslant$ then the
  semigroup $\mathscr{B}^+(G)$ is isomorphic to the bicyclic  monoid ${\mathscr{C}}(p,q)$ and the semigroup $\mathscr{B}(G)$ is isomorphic to the extended bicyclic semigroup $\mathscr{C}_{\mathbb{Z}}$ (see \cite{Fihel-Gutik-2011});

  \item[$(2)$] if $G$ is the additive group of real numbers
  $(\mathbb{R},+)$ with usual linear order $\leqslant$ then the
  semigroup $\mathscr{B}(G)$ is isomorphic to
  $B^{2}_{(-\infty,\infty)}$ (see \cite{Korkmaz-1997, Korkmaz-2009})
  and the semigroup $\mathscr{B}^+(G)$ is isomorphic to
  $B^{1}_{[0,\infty)}$ (see
  \cite{Ahre-1981, Ahre-1983, Ahre-1984, Ahre-1986, Ahre-1989}), \; and

  \item[$(3)$] the semigroup $\mathscr{B}(G)$ is isomorphic to
  the semigroup $S(G)$ which is defined in \cite{Fotedar-1974, Fotedar-1978}.
\end{itemize}

In the paper \cite{Gutik-Pagon-Pavlyk-2011} semigroups $\mathscr{B}(G)$ and $\mathscr{B}^+(G)$ are studied for a linearly ordered group $G$. That paper describes Green's relations on $\mathscr{B}(G)$ and $\mathscr{B}^+(G)$ and their bands, and shows that these are bisimple. Also in \cite{Gutik-Pagon-Pavlyk-2011} it is proved that for a commutative linearly ordered group $G$ all non-trivial congruences on the semigroups $\mathscr{B}(G)$ and $\mathscr{B}^+(G)$ are group congruences if and only if the group $G$ is Archimedean; and the structure of group congruences on the semigroups $\mathscr{B}(G)$ and $\mathscr{B}^+(G)$ is described.

In this paper we present more general construction than the semigroups $\mathscr{B}(G)$ and $\mathscr{B}^+(G)$.  Namely, for
a linearly ordered group $G$ let us define a subset $A\subseteq G$ to be a \emph{shift-set} if for any
$x,y,z\in A$ with $y < x$ we get $x\cdot y^{-1}\cdot z\in A$. For any shift-set $A\subseteq G$ let
\begin{equation*}
  \mathscr{B}(A)=\left\{\alpha_b^a\colon G^+(a)\to G^+(b)\colon a,b\in A\right\}
\end{equation*}
be the semigroup of partial bijections defined by the formula
\begin{equation*}
    (x)\alpha_b^a=x\cdot a^{-1}\cdot b, \qquad \hbox{ for } \; x\in
    G^{+}(a).
\end{equation*}
The semigroup $\mathscr{B}(A)$ is isomorphic to the semigroup $S_A=A\times A$ endowed with the binary operation defined by  formula \eqref{formula-1.2}. For $A=G$ the semigroup $\mathscr{B}(A)$ coincides with $\mathscr{B}(G)$ and for $A=G^+$ it coincides with the semigroup $\mathscr{B}^+(G)$.

Later in this paper for a non-empty shift-set $A\subseteq G$ we identify the semigroup $\mathscr{B}(A)$  with the semigroup $S_A$  endowed with the multiplication
defined by formula \eqref{formula-1.2}. We observe that $\mathscr{B}(A)$ is an inverse subsemigroup of $\mathscr{B}(G)$ for any non-empty shift-set $A$ of a linearly ordered group $G$. Moreover, the results of \cite{Gutik-Pagon-Pavlyk-2011} imply that an element $(a,b)$ of $\mathscr{B}(A)$ is an idempotent iff $a=b$, and $(b,a)$ is inverse of $(a,b)$ in $\mathscr{B}(G)$.

We recall that a topological space $X$ is said to be
\begin{itemize}
  \item \emph{locally compact}, if every point $x\in X$ has an open neighbourhood with the compact closure;
  \item \emph{\v{C}ech-complete}, if $X$ is Tychonoff and $X$ is a $G_\delta$-set in its the \v{C}ech-Stone compactification;
  \item \emph{Baire}, if for each sequence $A_1, A_2,\ldots, A_i,\ldots$ of open dense subsets of $X$ the intersection $\bigcap_{i=1}^\infty A_i$ is dense in $X$.
\end{itemize}
Every Hausdorff locally compact space is \v{C}ech-complete, and every \v{C}ech-complete space is Baire (see \cite{Engelking-1989}).

\textcolor[rgb]{0.00,0.00,0.00}{A \emph{semitopological} (\emph{topological}) \emph{semigroup} is a topological space with a separately continuous (jointly continuous) semigroup operation.}

A topology $\tau$ on a semigroup $S$ is called:
\begin{itemize}
  \item \emph{semigroup} if $(S,\tau)$ is a topological semigroup;
  \item \emph{shift-continuous} if $(S,\tau)$ is a semitopological semigroup.
\end{itemize}

The bicyclic monoid admits only the discrete semigroup Hausdorff topology and if a topological semigroup $S$ contains it as a dense subsemigroup then ${\mathscr{C}}(p,q)$ is an open subset of $S$~\cite{Eberhart-Selden-1969}. We observe that the openness of ${\mathscr{C}}(p,q)$ in its closure easily follows from the non-topologizability of the bicyclic monoid, because the discrete subspace $D$ is open in its closure $\overline{D}$ in any $T_1$-space containing $D$. Bertman and  West in \cite{Bertman-West-1976} extended this result for the case of Hausdorff semitopological semigroups. Stable and $\Gamma$-compact topological semigroups do not contain the bicyclic monoid~\cite{Anderson-Hunter-Koch-1965, Hildebrant-Koch-1986}. The problem of embedding the bicyclic monoid into compact-like topological semigroups was studied in \cite{Banakh-Dimitrova-Gutik-2009, Banakh-Dimitrova-Gutik-2010, Gutik-Repovs-2007}.
Independently Taimanov in \cite{Taimanov-1973} constructed a semigroup $\mathfrak{A}_\kappa$ of cardinality $\kappa$ which admits only the discrete semigroup topology. Also, Taimanov \cite{Taimanov-1975} gave sufficient conditions on a commutative semigroup to have a non-discrete semigroup topology. In the paper \cite{Gutik-2016} it was showed that for the Taimanov semigroup $\mathfrak{A}_\kappa$ from \cite{Taimanov-1973} the following conditions hold:
every $T_1$-topology $\tau$ on the semigroup $\mathfrak{A}_\kappa$ such that $(\mathfrak{A}_\kappa,\tau)$  is a topological semigroup is discrete; for every $T_1$-topological semigroup which contains $\mathfrak{A}_\kappa$ as a subsemigroup, $\mathfrak{A}_\kappa$ is a closed subsemigroup of $S$;
and every homomorphic non-isomorphic image of $\mathfrak{A}_\kappa$ is a zero-semigroup.
Also, in the paper \cite{Fihel-Gutik-2011} it is proved that the discrete topology is the unique shift-continuous Hausdorff topology on the extended bicyclic semigroup $\mathscr{C}_{\mathbb{Z}}$. Also, for many ($0$-) bisimple semigroups of transformations $S$ the following statement holds: \emph{every shift-continuous Hausdorff Baire (in particular locally compact) topology on $S$ is discrete} (see \cite{Chuchman-Gutik-2010, Chuchman-Gutik-2011, Gutik-Pozdnyakova-2014, Gutik-Repovs-2011, Gutik-Repovs-2012}). In the paper \cite{Mesyan-Mitchell-Morayne-Peresse-2016} Mesyan, Mitchell, Morayne and P\'{e}resse showed that if $E$ is a finite graph, then the only locally compact Hausdorff semigroup topology on the graph inverse semigroup $G(E)$ is the discrete topology. In \cite{Bardyla-Gutik-2016} it was proved that the conclusion of this statement also holds for graphs $E$ consisting of one vertex and infinitely many loops (i.e., infinitely generated polycyclic monoids). Amazing dichotomy for the bicyclic monoid with adjoined zero $\mathscr{C}^0={\mathscr{C}}(p,q)\sqcup\{0\}$ was proved in \cite{Gutik-2015}: every Hausdorff locally compact semitopological bicyclic monoid $\mathscr{C}^0$ with adjoined zero is either compact or discrete.
The above dichotomy was extended by Bardyla in \cite{Bardyla-2016} to locally compact $\lambda$-polycyclic semitopological monoids and to locally compact semitopological interassociates of the bicyclic monoid \cite{Gutik-Maksymyk-2016}.

For a linearly ordered group $G$ and a non-empty shift-set $A$ of $G$, the natural partial order and solutions of equations on the semigroup $\mathscr{B}(A)$ are described. We study topologizations of the semigroup $\mathscr{B}(A)$. In particular, we show that for an arbitrary countable linearly ordered group $G$ and a non-empty shift-set $A$ of $G$,  every Baire shift-continuous $T_1$-topology $\tau$ on $\mathscr{B}(A)$ is discrete. Also we prove that for an arbitrary linearly non-densely ordered group $G$ and a non-empty shift-set $A$ of $G$, every shift-continuous Hausdorff topology $\tau$ on the semigroup $\mathscr{B}(A)$ is discrete, and hence $\left(\mathscr{B}(A),\tau\right)$ is a discrete subspace of any Hausdorff semitopological semigroup which contains $\mathscr{B}(A)$ as a subsemigroup.


\section{Solutions of some equations and the natural partial order on the semigroup $\mathscr{B}(A)$}

It is well known that every inverse semigroup $S$ admits the \emph{natural partial order}:
\begin{equation*}
    s\preccurlyeq t \quad \hbox{~if and only if~} \quad s=et \quad \hbox{for some} \quad e\in E(S).
\end{equation*}
This order induces the natural partial order on the semilattice $E(S)$, and for arbitrary $s,t\in S$ the following conditions are equivalent:
\begin{equation}\label{formula-2.1}
    (\alpha)~s\preccurlyeq t; \qquad (\beta)~s=ss^{-1}t; \qquad (\gamma)~s=ts^{-1}s,
\end{equation}
(see \cite[Chapter~3]{Lawson-1998}.)

\begin{proposition}\label{proposition-2.1-0} Let $G$ be a linearly
ordered group and $A$ be a non-empty shift-set in $G$. Then the following assertions hold:
\begin{itemize}
  \item[$(i)$] if $(g,g),(h,h)\in E(\mathscr{B}(A))$        then $(g,g)\preccurlyeq(h,h)$ if and only if
       $g\geqslant h$ in $A$;

  \item[$(ii)$] the semilattice $E(\mathscr{B}(A))$ is isomorphic to $A$, considered as a $\vee$-semilattice under the isomorphisms
       $\mathfrak{i}\colon E(\mathscr{B}(A))\to A$, $\mathfrak{i}\colon(g,g)\mapsto g$;

  \item[$(iii)$] $(g,h)\mathscr{R}(k,l)$ in $\mathscr{B}(A)$  if and only if $g=k$ in $A$;

  \item[$(iv)$] $(g,h)\mathscr{L}(k,l)$ in $\mathscr{B}(A)$  if and only if $h=l$ in $A$;

  \item[$(v)$] $(g,h)\mathscr{H}(k,l)$ in $\mathscr{B}(A)$  if and only if $g=k$ and $h=l$ in $A$, and
       hence every $\mathscr{H}$-class in $\mathscr{B}(A)$ is a singleton;

  \item[$(vi)$] $\mathscr{B}(A)$ is a bisimple semigroup and hence it is simple.
\end{itemize}
\end{proposition}

\begin{proof}
Assertions $(i)$ and $(ii)$ are trivial, $(iii)-(v)$ follow  from Proposition~2.1 of \cite{Gutik-Pagon-Pavlyk-2011} and Proposition~3.2.11 of \cite{Lawson-1998} and $(vi)$ follows from Proposition~3.2.5 of \cite{Lawson-1998}.
\end{proof}

Later we need the following lemma, which describes the natural partial order on the semigroup $\mathscr{B}(A)$:

\begin{lemma}\label{lemma-2.1}
Let $G$ be a linearly ordered group and $A$ be a non-empty shift-set in $G$. Then for arbitrary elements $(a,b),(c,d)\in\mathscr{B}(A)$ the following conditions are equivalent:
\begin{itemize}
  \item[$(i)$] $(a,b)\preccurlyeq(c,d)$ in $\mathscr{B}(A)$;
  \item[$(ii)$] $a^{-1}\cdot b=c^{-1}\cdot d$ and $a\geqslant c$ in $A$;
  \item[$(iii)$] $b^{-1}\cdot a=d^{-1}\cdot b$ and $b\geqslant d$ in $A$.
\end{itemize}
\end{lemma}

\begin{proof}
$(i)\Rightarrow(ii)$ The equivalence of conditions $(\alpha)$ and $(\beta)$ in \eqref{formula-2.1} implies that $(a,b)\preccurlyeq(c,d)$ in $\mathscr{B}(A)$ if and only if $(a,b)=(a,b)(a,b)^{-1}(c,d)$. Therefore we have that
\begin{equation*}
    (a,b)=(a,b)(a,b)^{-1}(c,d)=(a,b)(b,a)(c,d)=(a,a)(c,d)=
    \left\{
      \begin{array}{ll}
        (c\cdot a^{-1}\cdot a,d), & \hbox{if~}\; a<c;\\
        (c,d),                    & \hbox{if~}\; a=c;\\
        (a,a\cdot c^{-1}\cdot d), & \hbox{if~}\; a>c.
      \end{array}
    \right.
\end{equation*}
This implies that
\begin{equation*}
    (a,b)=
    \left\{
      \begin{array}{ll}
        (c,d),                    & \hbox{if~}\; a<c;\\
        (c,d),                    & \hbox{if~}\; a=c;\\
        (a,a\cdot c^{-1}\cdot d), & \hbox{if~}\; a>c,
      \end{array}
    \right.
\end{equation*}
and hence the condition $(a,b)\preccurlyeq(c,d)$ in $\mathscr{B}(A)$ implies that $a^{-1}\cdot b=c^{-1}\cdot d$ and $a\geqslant c$ in $A$.

$(ii)\Rightarrow(i)$ Fix arbitrary $(a,b),(c,d)\in\mathscr{B}(A)$ such that $a^{-1}\cdot b=c^{-1}\cdot d$ and $a\geqslant c$ in $A$. Then we have that
\begin{equation*}
    (a,b)(a,b)^{-1}(c,d)=(a,b)(b,a)(c,d)=(a,a)(c,d)=(a,a\cdot c^{-1}\cdot d)=(a,b),
\end{equation*}
and hence $(a,b)\preccurlyeq(c,d)$ in $\mathscr{B}(A)$.

The proof of the equivalence $(ii)\Leftrightarrow(iii)$ is trivial.
\end{proof}

The definition of the semigroup operation in $\mathscr{B}(A)$ implies that $(a,b)=(a,c)(c,d)(d,b)$ for arbitrary elements $a,b,c,d$ of $A$. The following two propositions give amazing descriptions of solutions of some equations in the semigroup $\mathscr{B}(A)$.

\begin{proposition}\label{proposition-2.2}
Let $G$ be a linearly ordered group, $A$ be a non-empty shift-set in $G$ and $a,b,c,d$ be arbitrary elements of $A$. Then the following con\-di\-tions hold:
\begin{itemize}
  \item[$(i)$] $(a,b)=(a,c)(x,y)$ for $(x,y)\in\mathscr{B}(A)$ if and only if $(c,b)\preccurlyeq(x,y)$ in $\mathscr{B}(A)$;
  \item[$(ii)$] $(a,b)=(x,y)(d,b)$ for $(x,y)\in\mathscr{B}(A)$ if and only if $(a,d)\preccurlyeq(x,y)$ in $\mathscr{B}(A)$;
  \item[$(iii)$] $(a,b)=(a,c)(x,y)(d,b)$ for $(x,y)\in\mathscr{B}(A)$ if and only if $(c,d)\preccurlyeq(x,y)$ in $\mathscr{B}(A)$.
\end{itemize}
\end{proposition}

\begin{proof}
$(i)$~$(\Rightarrow)$
Suppose that $(a,b)=(a,c)(x,y)$ for some $(x,y)\in\mathscr{B}(A)$. Then we have that
\begin{equation*}
    (a,c)(x,y)=
\left\{
      \begin{array}{ll}
        (a,c\cdot x^{-1}\cdot y), & \hbox{if~}\; c>x;\\
        (a,y),                    & \hbox{if~}\; c=x;\\
        (x\cdot c^{-1}\cdot a,y), & \hbox{if~}\; c<x.
      \end{array}
\right.
\end{equation*}
Then in the case when $c>x$ we get that $b=c\cdot x^{-1}\cdot y$ and hence Lemma~\ref{lemma-2.1} implies that $(c,b)\preccurlyeq(x,y)$ in $\mathscr{B}(A)$. Also, in the case when $c=x$ we have that $b=y$, which implies the inequality $(c,b)\preccurlyeq(x,y)$ in $\mathscr{B}(A)$. The case $c<x$ does not hold because the group operation on $G$ implies that $x\cdot c^{-1}\cdot a<a$.

$(\Leftarrow)$
Suppose that the relation $(c,b)\preccurlyeq(x,y)$ holds in $\mathscr{B}(A)$. Then by Lemma~\ref{lemma-2.1} we have that $c^{-1}\cdot b=x^{-1}\cdot y$ and $c\geqslant x$ in $A$, and hence the semigroup operation of $\mathscr{B}(A)$ implies that
\begin{equation*}
    (a,c)(x,y)=(a,c\cdot x^{-1}\cdot y)=(a,c\cdot c^{-1}\cdot b)=(a,b).
\end{equation*}

The proof of statement $(ii)$ is similar to statement $(i)$.

$(iii)$~$(\Rightarrow)$
Suppose that $(a,b)=(a,c)(x,y)(d,b)$ for some $(x,y)\in\mathscr{B}(A)$. Then we have that
\begin{equation*}
    (a,c)(x,y)=
\left\{
      \begin{array}{ll}
        (a,c\cdot x^{-1}\cdot y), & \hbox{if~}\; c>x;\\
        (a,y),                    & \hbox{if~}\; c=x;\\
        (x\cdot c^{-1}\cdot a,y), & \hbox{if~}\; c<x.
      \end{array}
\right.
\end{equation*}
Therefore,
\begin{itemize}
  \item[(a)] if $c>x$ then
  \begin{equation*}
    (a,c)(x,y)(d,b)=(a,c\cdot x^{-1}\cdot y)(d,b)=
  \left\{
      \begin{array}{ll}
        (a,c\cdot x^{-1}\cdot y\cdot d^{-1}\cdot b), & \hbox{if~}\; c\cdot x^{-1}\cdot y>d;\\
        (a,b),                                       & \hbox{if~}\; c\cdot x^{-1}\cdot y=d;\\
        (d\cdot y^{-1}\cdot x\cdot c^{-1}\cdot a,b), & \hbox{if~}\; c\cdot x^{-1}\cdot y<d,
      \end{array}
\right.
  \end{equation*}
  \item[(b)] if $c=x$ then
  \begin{equation*}
    (a,c)(x,y)(d,b)=(a,y)(d,b)=
  \left\{
      \begin{array}{ll}
        (a,y\cdot d^{-1}\cdot b), & \hbox{if~}\; y>d;\\
        (a,b),                    & \hbox{if~}\; y=d;\\
        (d\cdot y^{-1}\cdot a,b), & \hbox{if~}\; y<d,
      \end{array}
  \right.
  \end{equation*}
  \item[(c)] if $c<x$ then
  \begin{equation*}
    (a,c)(x,y)(d,b)=(x\cdot c^{-1}\cdot a,y)(d,b)=
  \left\{
      \begin{array}{ll}
        (x\cdot c^{-1}\cdot a,y\cdot d^{-1}\cdot b), & \hbox{if~}\; y>d;\\
        (x\cdot c^{-1}\cdot a,b),                    & \hbox{if~}\; y=d;\\
        (d\cdot y^{-1}\cdot x\cdot c^{-1}\cdot a,b), & \hbox{if~}\; y<d.
      \end{array}
  \right.
  \end{equation*}
\end{itemize}
Then the equality $(a,b)=(a,c)(x,y)(d,b)$ implies that
\begin{itemize}
  \item[(a)] if $c>x$ then $c\cdot x^{-1}\cdot y\cdot d^{-1}=e$ in $G$;

  \item[(b)] if $c=x$ then $y=d$;
\end{itemize}
and case (c) does not hold. Hence by Lemma~\ref{lemma-2.1} we get that $(c,d)\preccurlyeq(x,y)$ in $\mathscr{B}(A)$.

$(\Leftarrow)$
Suppose that the relation $(c,d)\preccurlyeq(x,y)$ holds in $\mathscr{B}(A)$. Then by Lemma~\ref{lemma-2.1} we have that $c^{-1}\cdot d=x^{-1}\cdot y$ and $c\geqslant x$ in $A$, and hence the semigroup operation of $\mathscr{B}(A)$ implies
\begin{equation*}
\begin{split}
  (a,c)(x,y)(d,b)= &\; (a,c)(x,y)(c\cdot x^{-1}\cdot y,b)=\\
    = &\; (a,c)(c\cdot x^{-1}\cdot y\cdot y^{-1}\cdot x,b)=\\
    = &\; (a,c)(c\cdot x^{-1}\cdot x,b)=\\
    = &\; (a,c)(c,b)=\\
    = &\; (a,b),
\end{split}
\end{equation*}
because $c\cdot x^{-1}\cdot y\geqslant y$ in $A$.
\end{proof}

\begin{proposition}\label{proposition-2.3}
Let $G$ be a linearly ordered group, $A$ be a non-empty shift-set in $G$ and $a,b,c,d$ be arbitrary elements of $A$. Then the following statements hold:
\begin{itemize}
  \item[$(i)$] if $a<c$ in $G$ then the equation $(a,b)=(c,d)(x,y)$ has no solutions in $\mathscr{B}(A)$;
  \item[$(ii)$] if $a>c$ in $A$ then the equation $(a,b)=(c,d)(x,y)$ has the unique solution $(x,y)=(a\cdot c^{-1}\cdot d,b)$ in $\mathscr{B}(A)$;
  \item[$(iii)$] $(a,b)=(a,d)(x,y)$ has the unique solution $(x,y)=(d,b)$ in $\mathscr{B}(A)$;
  \item[$(iv)$] if $b<d$ in $A$ then the equation $(a,b)=(x,y)(c,d)$ has no solutions in $\mathscr{B}(A)$;
  \item[$(v)$] if $b>d$ in $A$ then the equation $(a,b)=(x,y)(c,d)$ has the unique solution $(x,y)=(a,b\cdot d^{-1}\cdot c)$ in $\mathscr{B}(A)$;
  \item[$(v)$] the equation $(a,b)=(x,y)(c,b)$ has the unique solution $(x,y)=(a,c)$ in $\mathscr{B}(A)$.
\end{itemize}
\end{proposition}

\begin{proof}
$(i)$ Assume that $a<c$. Then formula \eqref{formula-1.2} implies that $d<x$ in $A$ and hence $(a,b)=(x\cdot d^{-1}\cdot c,y)$. This implies that $a=x\cdot d^{-1}\cdot c$ and $b=y$. Since $d<x$, the equality $a=x\cdot d^{-1}\cdot c$ implies $a>c$, which contradicts the assumption of statement $(i)$.

$(ii)$ Assume that $a>c$. Then formula \eqref{formula-1.2} implies that $d<x$ in $A$ and hence we have that $(a,b)=(x\cdot d^{-1}\cdot c,y)$. This implies the equalities $x=a\cdot c^{-1}\cdot d$ and $y=b$.

$(iii)$ follows from formula \eqref{formula-1.2}.

The proofs of statements $(iv)$, $(v)$ and $(vi)$ are dual to the proofs of $(i)$, $(ii)$ and $(iii)$, respectively.
\end{proof}

Later we need the following proposition which follows from formula \eqref{formula-1.2} and describes right and left principal ideals in the semigroup  $\mathscr{B}(A)$ for a non-empty shift-set $A$ in $G$.

\begin{proposition}\label{proposition-2.7}
Let $G$ be a linearly ordered group and $A$ be a non-empty shift-set in $G$. Then the following con\-di\-ti\-ons hold:
\begin{itemize}
  \item[$(i)$] $(a,a)\mathscr{B}(A)=\left\{(x,y)\in\mathscr{B}(A)\colon x\geqslant a \hbox{~in~} A\right\}$;
  \item[$(ii)$] $\mathscr{B}(A)(a,a)=\left\{(x,y)\in\mathscr{B}(A)\colon y\geqslant a \hbox{~in~} A\right\}$.
\end{itemize}
\end{proposition}


\section{On topologizations of the semigroup $\mathscr{B}(A)$}

It is obvious that every left (right) topological group $G$ with an isolated point is discrete. This implies that every countable $T_1$-Baire left (right) topological group is a discrete space, too. Later we shall show that the similar statement holds for a Baire semitopological semigroup $\mathscr{B}(A)$ over a non-empty shift-set $A$ of a countable linearly ordered group  $G$.

For an arbitrary element $(a,b)$ of the semigroup $\mathscr{B}(A)$ we denote
\begin{equation*}
{\uparrow}_{\preccurlyeq}(a,b)=\left\{(x,y)\in\mathscr{B}(A) \colon (a,b)\preccurlyeq(x,y)\right\}.
\end{equation*}

\begin{lemma}\label{lemma-3.1}
Let $G$ be a linearly ordered group, $A$ be a non-empty shift-set in $G$, and $\tau$ be a shift-continuous topology on $\mathscr{B}(A)$ such that $\left(\mathscr{B}(A),\tau\right)$ contains an isolated point. Then the space $\left(\mathscr{B}(A),\tau\right)$ is discrete.
\end{lemma}

\begin{proof}
Suppose that $(a,b)$ is an isolated point of the topological space $\left(\mathscr{B}(A),\tau\right)$. Assume that for an arbitrary $u\in A$ there exists $c\in A$ such that $u>c$. Since $A$ is a shift-set, $d=c\cdot u^{-1}\cdot b<b$ in $A$. By Proposition~\ref{proposition-2.3}$(v)$, the equation $(a,b)=(x,y)(c,d)$ has the unique solution
\begin{equation*}
(x,y)=\left(a,b\cdot d^{-1}\cdot c\right)=\left(a,b\cdot (c\cdot u^{-1}\cdot b)^{-1}\cdot c\right)=\left(a,b\cdot b^{-1}\cdot u\cdot c^{-1}\cdot c\right)=(a,u)
\end{equation*}
in $\mathscr{B}(A)$. If $u$ is the smallest element of $A$, then By Proposition~\ref{proposition-2.3}$(vi)$ the equation $(a,b)=(x,y)(u,b)$ has the unique solution $(x,y)=(a,u)$. In both cases the continuity of right translations in $\left(\mathscr{B}(A),\tau\right)$ implies that for
arbitrary $u\in A$ the pair  $(a,u)$ is an isolated point of the topological space $\left(\mathscr{B}(A),\tau\right)$ for arbitrary $u\in A$.

Fix an arbitrary element $v$ of $A$. Assume that there exists $d\in A$ such that $d<v$. Since $A$ is a shift-set, $c=d\cdot v^{-1}\cdot a<a$ in $A$. Then by Proposition~\ref{proposition-2.3}$(ii)$,  the equation $(a,u)=(c,d)(x,y)$ has the unique solution
\begin{equation*}
(x,y)=\left(a\cdot c^{-1}\cdot d,u\right)=\left(a\cdot (d\cdot v^{-1}\cdot a)^{-1}\cdot d,u\right)=\left(a\cdot a^{-1}\cdot v\cdot d^{-1}\cdot d,u\right)=(v,u)
\end{equation*}
in $\mathscr{B}(A)$. If $v$ is the smallest element of $A$, then by Proposition~\ref{proposition-2.3}$(iii)$, the equation $(a,u)=(a,v)(x,y)$ has the unique solution $(x,y)=(v,u)$. Since $(a,u)$ is an isolated point of $\left(\mathscr{B}(G),\tau\right)$, in both cases the continuity of left translations in $\mathscr{B}(G)$ implies that for arbitrary $u\in A$ the pair  $(v,u)$ is an isolated point of the topological space $\left(\mathscr{B}(G),\tau\right)$ for arbitrary $u\in G$. This completes the proof of the lemma.
\end{proof}

\begin{theorem}\label{theorem-3.2}
Let $A$ be a countable non-empty shift-set in a linearly ordered group $G$,  and $\tau$ be a $T_1$-Baire shift-continuous topology on $\mathscr{B}(A)$. Then the topological space $\left(\mathscr{B}(A),\tau\right)$ is discrete.
\end{theorem}

\begin{proof}
By Proposition~1.30 of \cite{Haworth-McCoy-1977} every countable Baire $T_1$-space contains a dense subspace of isolated points, and hence the space $\left(\mathscr{B}(A),\tau\right)$ contains an isolated point. Then we apply Lemma~\ref{lemma-3.1}.
\end{proof}

Theorem~\ref{theorem-3.2} implies the following corollary:

\begin{corollary}\label{corollary-3.3}
Let $A$ be a countable non-empty shift-set in a linearly ordered group $G$, and $\tau$ be a shift-continuous \v{C}ech complete $($or locally compact$)$ $T_1$-topology on $\mathscr{B}(A)$. Then the topological space $\left(\mathscr{B}(A),\tau\right)$ is discrete.
\end{corollary}

\begin{remark}\label{remark-3.4}
Let $\mathbb{R}$ be the set of reals with usual topology. It is obvious that $S_{\mathbb{R}}=\mathbb{R}\times\mathbb{R}$ with the semigroup operation
  \begin{equation*}\label{formula-3.1}
  (a,b)\cdot(c,d)=
  \left\{
    \begin{array}{ll}
      (a-b+c,d), & \hbox{if }  \; b<c;\\
      (a,d),     & \hbox{if } \; b=c; \\
      (a,b-c+d), & \hbox{if } \; b>c,
    \end{array}
  \right.
  \end{equation*}
is isomorphic to the semigroup $\mathscr{B}(\mathbb{R})$, where $\mathbb{R}$ is the additive group of reals with usual linear order. Then simple verifications show that $S$ with the product topology $\tau_p$ is a topological inverse semigroup (also, see \cite{Korkmaz-1997, Korkmaz-2009}). Then the subspace $S_{\mathbb{Q}}=\{(x,y)\in S_{\mathbb{R}}\colon x \hbox{~and~} y \hbox{~are rational}\}$ with the induced semigroup operation from $S$ is a countable non-discrete non-Baire topological inverse subsemigroup of $(S,\tau_p)$. Also, the same we get in the case of subsemigroup $S^+_{\mathbb{Q}}=\{(x,y)\in S_{\mathbb{Q}}\colon x\geqslant0 \hbox{~and~} y\geqslant0 \}$ of $(S,\tau_p)$ (see \cite{Ahre-1981, Ahre-1983, Ahre-1984, Ahre-1986, Ahre-1989}). The above arguments show that the condition in Theorem~\ref{theorem-3.2} that $\tau$ is a $T_1$-Baire topology is essential.
\end{remark}

Recall that a linearly ordered group $G$ is said to be \emph{densely ordered} if for every positive element $g\in G$ there exists a positive element $h\in G$ such that $h<g$.

\begin{remark}\label{remark-3.5}
It is obviously that for a linearly ordered group $G$ the following conditions are equivalent:
\begin{enumerate}
  \item[$(i)$] $G$ is not densely ordered;
  \item[$(ii)$] for every $g\in G$ there exists a unique $g^+\in G$ such that $G^+(g)\setminus G^+(g^+)=\{g\}$;
  \item[$(iii)$] for every $g\in G$ there exists a unique $g^-\in G$ such that $G^-(g)\setminus G^-(g^-)=\{g\}$, where $G^-(g)$ is the \emph{negative cone on} the element $g$, i.e., $G^-(g)=\{x\in G\colon x\leqslant g\}$.
\end{enumerate}
\end{remark}

Later for a linearly ordered group $G$ which is not densely ordered and an arbitrary element $g$ of a non-empty shift-set $A$ in $G$ by $g^+$ (resp., $g^-$) we denote the minimum (resp., maximum) element of the set $G^+(g)\setminus\{g\}\cap A$ (resp., $G^-(g)\setminus\{g\}\cap A$).

\begin{theorem}\label{theorem-3.6}
Let $G$ be a linearly ordered group which is not densely ordered and $A$ be a non-empty shift-set in $G$.
Then every shift-continuous Hausdorff topology $\tau$ on the semigroup
$\mathscr{B}(A)$ is discrete, and hence $\mathscr{B}(A)$ is a discrete
subspace of any semitopological semigroup which contains
$\mathscr{B}(A)$ as a subsemigroup.
\end{theorem}

\begin{proof}
We fix an arbitrary idempotent $(a,a)$ of the semigroup
$\mathscr{B}(A)$ and suppose that $(a,a)$ is a
non-isolated point of the topological space
$(\mathscr{B}(A),\tau)$. Since the maps
$\lambda_{(a,a)}\colon
\mathscr{B}(A)\rightarrow\mathscr{B}(A)$ and
$\rho_{(a,a)}\colon \mathscr{B}(A)\rightarrow
\mathscr{B}(A)$ defined by the formulae
$\left((x,y)\right)\lambda_{(a,a)}=(a,a)(x,y)$ and
$\left((x,y)\right)\rho_{(a,a)}=(x,y)(a,a)$ are continuous
retractions, we conclude that $(a,a)\mathscr{B}(A)$ and
$\mathscr{B}(A)(a,a)$ are closed subsets in the
topological space $(\mathscr{B}(A),\tau)$ (see \cite[Exercise~1.5.C]{Engelking-1989}). For an arbitrary element $b$ of the shift-set $A$ in the linearly ordered group $G$ we put
\begin{equation*}
    \textsf{DL}_{(b,b)}\left[(b,b)\right]= \left\{(x,y)\in\mathscr{B}(A) \colon (x,y)(b,b)=(b,b)\right\}.
\end{equation*}
Lemma~\ref{lemma-2.1} and Proposition~\ref{proposition-2.2} imply that
\begin{equation*}
    \textsf{DL}_{(b,b)}\left[(b,b)\right]={\uparrow}_{\preccurlyeq}(b,b)= \left\{(x,x)\in\mathscr{B}(A) \colon x\leqslant b \; \hbox{~in~} \; A\right\},
\end{equation*}
and since right translations are continuous maps in $(\mathscr{B}(A),\tau)$ we get that $\textsf{DL}_{(b,b)}\left[(b,b)\right]$ is a closed subset of the
topological space $(\mathscr{B}(A),\tau)$ for every $b\in A$. Then there exists an open neighbourhood $W_{(a,a)}$ of the point $(a,a)$ in the
topological space $(\mathscr{B}(A),\tau)$ such that
\begin{equation*}
    W_{(a,a)}\subseteq\mathscr{B}(A)\setminus \big((a^+,a^+)\mathscr{B}(A)\cup \mathscr{B}(A)(a^+,a^+)\cup \textsf{DL}(a^-,a^-)\big).
\end{equation*}
Since $(\mathscr{B}(A),\tau)$ is a semitopological semigroup, we conclude that there exists an open neighbourhood $V_{(a,a)}$ of the idempotent $(a,a)$ in the topological space $(\mathscr{B}(A),\tau)$ such that the following conditions hold:
\begin{equation*}
    V_{(a,a)}\subseteq W_{(a,a)}, \qquad
    (a,a)\cdot V_{(a,a)}\subseteq W_{(a,a)} \qquad \hbox{and} \qquad
    V_{(a,a)}\cdot(a,a) \subseteq W_{(a,a)}.
\end{equation*}
Hence at least one of the following conditions holds:
\begin{itemize}
  \item[$(a)$] the neighbourhood $V_{(a,a)}$ contains infinitely many points $(x,y)\in\mathscr{B}(A)$ such that $x<y\leqslant a$ in the set $A$; \; or

  \item[$(b)$] the neighbourhood $V_{(a,a)}$ contains infinitely many points $(x,y)\in\mathscr{B}(A)$ such that $y<x\leqslant a$ in the set $A$.
\end{itemize}
In case $(a)$ by Proposition~\ref{proposition-2.2} we have that
\begin{equation*}
    (a,a) (x,y)=\left(a,a\cdot x^{-1}\cdot y\right)\notin W_{(a,a)}
\end{equation*}
because $x^{-1}\cdot y\geqslant e$ in $G$, and in case $(b)$ by Proposition~\ref{proposition-2.2} we have that
\begin{equation*}
    (x,y) (a,a)=\left(a\cdot y^{-1}\cdot x,a\right)\notin W_{(a,a)}
\end{equation*}
because $y^{-1}\cdot x\geqslant e$ in $G$ which contradicts the separate continuity of the semigroup operation in $(\mathscr{B}(A),\tau)$. The obtained
contradiction implies that the set $V_{(a,a)}$ is a singleton, and
hence the idempotent $(a,a)$ is an isolated point of the topological
space $(\mathscr{B}(A),\tau)$.

Now, we apply Lemma~\ref{lemma-3.1} and get that the topological space $(\mathscr{B}(A),\tau)$ is discrete.
\end{proof}

Theorem~\ref{theorem-3.6} implies the following three corollaries:

\begin{corollary}\label{corollary-3.7}
Let $G$ be a linearly ordered group which is not densely ordered and $A$ be a non-empty shift-set in $G$. Then every semigroup Hausdorff topology $\tau$ on the semigroup $\mathscr{B}(A)$  is discrete.
\end{corollary}

\begin{corollary}[\cite{Fihel-Gutik-2011}]\label{corollary-3.8}
Every shift-continuous Hausdorff topology $\tau$ on the bicyclic extended semigroup $\mathscr{C}_{\mathbb{Z}}$ is discrete.
\end{corollary}

\begin{corollary}[\cite{Bertman-West-1976, Eberhart-Selden-1969}]\label{corollary-3.9}
Every shift-continuous Hausdorff topology $\tau$ on the bicyclic monoid $\mathscr{C}(p,q)$ is discrete.
\end{corollary}

\section*{Acknowledgements}

The author acknowledges Taras Banakh and the referee for their important comments and suggestions.


\end{document}